# LARGE SCALE GEOMETRY OF 4-DIMENSIONAL CLOSED NONPOSITIVELY CURVED REAL ANALYTIC MANIFOLDS

MOHAMAD A. HINDAWI

ABSTRACT. We study the asymptotic cones of the universal covering spaces of closed 4-dimensional nonpositively curved real analytic manifolds. We show that the existence of nonstandard components in the Tits boundary, discovered by Christoph Hummel and Victor Schroeder [HS98], depends only on the quasi-isometry type of the fundamental group.

## 1. INTRODUCTION

The concept of the asymptotic cone was introduced by van den Dries and Wilkie [vdDW84] and by Gromov [Gro93]. It has been used by several authors to study the large scale geometry of certain spaces and to distinguish the quasi-isometry type of these spaces. Kleiner and Leeb [KL97b] used the asymptotic cone to prove the Margulis conjecture regarding the quasi-isometry rigidity of symmetric spaces of noncompact type. Kapovich and Leeb [KL95, KL97a] also used the asymptotic cone to study the quasi-isometry types of Haken 3-manifolds.

Nonpositively curved closed 4-dimensional real analytic manifolds form an interesting class of manifolds. Many examples of these manifolds have been constructed by Schroeder [Sch88, Sch89a] and later by Abresch and Schroeder [AS96]. The structure of flats in these manifolds was studied by Schroeder [Sch89b]. The Tits alternative of their fundamental groups was studied by Xiangdong Xie [Xie04].

Given a Hadamard space $X$, the Tits boundary $\partial_T X = (X(\infty), Td)$ is a metric space which reflects part of the asymptotic geometry of the space. The ideal boundary of a Hadamard space in general is *not* a quasi-isometry invariant. Croke and Kleiner gave the first example of such phenomenon. They constructed a pair of compact piecewise Euclidean 2-complexes with nonpositive curvature which are homeomorphic but whose universal covers have non-homeomorphic geometric boundaries [CK00]. This is in contrast with the case of Gromov hyperbolic

---

spaces, where quasi-isometric spaces have homeomorphic geometric boundaries. Buyalo [Buy98] and Croke and Kleiner [CK02] independently also showed that the Tits boundary is not equivariantly preserved under a quasi-isometry.

By investigating the Tits boundary of graph manifolds, Croke and Kleiner [CK02] discovered nontrivial connected components of the Tits boundary which are not subsets of unions of ideal boundaries of flats. These components are intervals of length less than $\pi$. Hummel and Schroeder [HS98] discovered similar components, which they called *nonstandard components,* in the Tits boundary of certain 4-dimensional nonpositively curved real analytic manifolds. It follows from [Sch89b] that the existence of these nonstandard components depends on the fundamental group of the 4-dimensional real analytic manifold.

In this paper we analyze, using some ideas of Kapovich and Leeb [KL95], the asymptotic cones of the universal cover of these manifolds, and show that the existence of nonstandard components depends only on the quasi-isometry type of the fundamental group of the manifold.

**Theorem 1.1.** *Let $X_1$ and $X_2$ be the universal covers of two closed 4-dimensional nonpositively curved real analytic manifolds. If $X_1$ and $X_2$ are quasi-isometric, then $\partial_T X_1$ contains a nonstandard component if and only if $\partial_T X_2$ does.*

If the Tits boundary of the universal covering space $X$ of a closed irreducible 4-dimensional nonpositively curved real analytic manifold does *not* contain any nonstandard components, we show in Section 4 that any asymptotic cone of $X$ is tree-graded in the sense of Druțu and Sapir [DS] (see Definition 1.10 in [DS]) with respect to a collection $\mathcal{F}$ of subsets, obtained as $\omega$-lim of higher rank submanifolds in $X$. The following corollary follows immediately from Theorem 8.5 of Druțu and Sapir [DS].

**Corollary 1.2.** *Let $X$ be the universal cover of a closed irreducible 4-dimensional nonpositively curved real analytic manifold $M$. If $\partial_T X$ does not contain any nonstandard components, then the fundamental group $\pi_1(M)$ is relatively hyperbolic with respect to the fundamental groups of the maximal higher rank submanifolds in $M$.*

The paper is organized as follows. Section 1 is an introduction. In Section 2 we recall the properties of 4-dimensional real analytic manifolds which we need through the rest of the paper. In Section 3 we show that for certain 4-dimensional real analytic manifolds all triangles are thin relative to a maximal higher rank submanifold. We use that to describe all flats in the asymptotic cone of these manifolds which will be done in Section 4. In Section 5 we give the proof of Theorem 1.1.

**Acknowledgment.** This paper is part of my thesis. I would like to thank my advisor Christopher B. Croke for many conversations related and unrelated to this paper and for his guiding and support over the years.

## 2. Background

Let $X^n$ be an $n$-dimensional Hadamard manifold, by that we mean a complete simply connected Riemannian manifold with nonpositive sectional curvature. By Cartan-Hadamard Theorem $X^n$ is diffeomorphic to $\mathbb{R}^n$. In fact the $\exp_p$ map at any point $p \in X$ is a diffeomorphism. For basic facts about Hadamard space, we refer the reader to [BGS85].

We denote the ideal boundary $X(\infty)$ of $X$ equipped with the cone topology by $\partial_\infty X$, and denote the ideal boundary with the Tits metric by $\partial_T X$. For any two different points $p, q \in X$, $\overline{pq}$, $\overrightarrow{pq}$ denote respectively the geodesic segment connecting $p$ to $q$, and the geodesic ray starting at $p$ and passing through $q$. By $\overrightarrow{pq}(\infty)$ we denote the limit point in $X(\infty)$ of the ray $\overrightarrow{pq}$.

Throughout the paper, with the exception of Proposition 3.3, $X$ will denote the universal cover of a closed irreducible 4-dimensional nonpositively curved real analytic manifold and we denote by $\Gamma$ its fundamental group.

Given a unit tangent vector $v$, the *rank* of $v$ is the dimension of the vector space of parallel Jacobi fields along the unique geodesic $c_v$ with $c'_v(0) = v$. If the isometry group of the Hadamard manifold X satisfies the duality condition, for example if X covers a manifold of finite volume, then parallel Jacobi fields can be integrated to produce flat strips, see Chapter IV in [Bal95] for details.

For any unit vector $v$, we write $v(\pm\infty) \in X(\infty)$ to denote the end points of the unique geodesic $c_v$. The parallel set of $c_v$ is denoted by $P_v$. This set consists of the union of all geodesics which are parallel to $c_v$. By analyticity, $P_v$ is a complete totally geodesic submanifold without boundary. Moreover $P_v$ splits isometrically as $Q \times \mathbb{R}$, where the $\mathbb{R}$-factor corresponds to the geodesic $c_v$.

A connected submanifold $F$ is said to be a *higher rank submanifold* of $X$ if it is totally geodesic submanifold with the property that every geodesic $c$ in $F$ has a parallel $c'$ in $F$ such that $c \neq c'$. We say $F$ is a *maximal higher rank submanifold* if it is not properly contained in any other higher rank submanifold. Given a unit vector $v$, $P_v$ is a higher rank submanifold unless $\text{rank}(v) = 1$. Schroeder [Sch89b] gave a complete description of the higher rank submanifolds in $X$. If $X$ is irreducible, then there are exactly three types, $\mathbb{R}^2$, $\mathbb{R}^3$, and $Q \times \mathbb{R}$, where $Q$ is a 2-dimensional visibility manifold. Moreover all the maximal higher rank submanifolds are closed under the action of the fundamental group, i.e., for any maximal higher rank submanifold $F$, $\Gamma_F = \{\phi \in \Gamma \mid \phi F = F\}$ operates with compact quotient on $F$.

We denote by $\mathcal{V}$ the set of all maximal higher rank submanifolds and by $\mathcal{W}$ the set of all maximal higher rank submanifolds of the form $Q \times \mathbb{R}$, where $Q$ is a 2-dimensional visibility manifold.

Schroeder [Sch89b] proved that modulo $\Gamma$ there are only finite number of maximal higher rank submanifolds. He also gave a description of the possible intersections of the maximal higher rank submanifolds. We now recall them. Given two different maximal higher rank submanifolds $F_1$ and $F_2$.

1. If $F_1 \approx \mathbb{R}^3$, then $F_1 \cap F_2 = \emptyset$.
2. If $F_1 \approx \mathbb{R}^2$ and $F_1 \cap F_2 \neq \emptyset$, then $F_2 \approx \mathbb{R}^2$ and $F_1 \cap F_2$ is a single point.
3. If $F_1 \in \mathcal{W}$ and $F_1 \cap F_2 \neq \emptyset$, then $F_2 \in \mathcal{W}$ and $F_1 \cap F_2$ is a 2-flat.

For any maximal higher rank submanifold $F \in \mathcal{V}$, $\partial_T F$ is a connected subset of $\partial_T X$. If $F \notin \mathcal{W}$, then $\partial_T F$ is isometric either to $S^1$ or $S^2$. If $F \in \mathcal{W}$, then $\partial_T F$ is a graph with two vertices and uncountable number of edges, where the length of each edge is $\pi$. The union of any two edges is a closed geodesic in $\partial_T F$ and it is the ideal boundary of a 2-flat in $F$. The two vertices are called singular points and they are precisely the end points of any singular geodesic, i.e., a geodesic of the form $\{q\} \times \mathbb{R} \subset Q \times \mathbb{R} = F$.

We now describe the possible intersections of the ideal boundaries of the maximal higher rank submanifolds. The only possible intersection is when $F_1, F_2 \in \mathcal{W}$ and $F_1 \cap F_2 = K$, where K is a 2-flat. In such case $\partial_T F_1 \cap \partial_T F_2 = \partial_T K = S^1$.

**Definition 2.1.** A connected component of $\partial_T X$ is called *standard* if it contains a boundary point of a flat and *nonstandard* if it is not a single point and not standard.

The main result of Hummel and Schroeder in [HS98] was to show that $\partial_T X$ contains nonstandard components precisely when there exist two maximal higher rank submanifolds $F_1, F_2 \in \mathcal{W}$ which intersect in a 2-flat, and that all the nonstandard components are intervals of length less than $\pi$.

## 3. Fat Triangles

Throughout this section, with the exception of Proposition 3.3, $X$ will denote the universal cover of a closed irreducible 4-dimensional nonpositively curved real analytic manifold which has no nonstandard component in $\partial_T X$.

If $X$ contains any higher rank submanifold, then it is not Gromov $\delta$-hyperbolic space for any $\delta \geq 0$. Easy examples of triangles which are not $\delta$-thin are triangles which lie in a flat in a higher rank submanifold. The goal of this section is to show that all "fat" triangles have to lie near a maximal higher rank submanifold, see Theorem 3.2. We begin with the following definition.

**Definition 3.1.** Let $\delta > 0$ and let $F$ be a maximal higher rank submanifold in $X$, a triangle $\Delta(p,q,r)$ is called $\delta$-thin relative to $F$, if every side is contained in a $\delta$-neighborhood of the union of the other two sides and $F$.

This definition resembles the definition of thin triangles relative to flats used by Hruska in [Hru04].

**Theorem 3.2.** *There exists some constant $\delta > 0$ such that any triangle in $X$ is $\delta$-thin relative to some maximal higher rank submanifold.*

Most of the proofs in this section depends on the following observation.

**Proposition 3.3.** *Let $X^n$ be an $n$-dimensional Hadamard manifold, not necessarily analytic. Let $p \in X$ and let $\overline{x_n y_n}$ be a sequence of geodesic segments such that $x_n, y_n$ converge respectively to $\xi_x, \xi_y \in X(\infty)$. If $d(p, \overline{x_n y_n})$ goes to infinity then $Td(\xi_x, \xi_y) \leq \pi$.*

*Proof.* The proof is easy and resembles the proof of Lemma 2.1 in [Bal82]. Assume that the proposition is not true. If $Td(\xi_x, \xi_y) > \pi$ then there exists a geodesic $\sigma$ in $X$ such that $\sigma(\infty) = \xi_x$ and $\sigma(-\infty) = \xi_y$. Notice that for any two points $p, q \in X$ the metrics $\angle_p$ and $\angle_q$ define the same topology on $\overline{X} = X \cup \partial_\infty X$. So without loss of generality and to simplify the notation, we assume that $p = \sigma(0)$.

Let $\sigma_n$ be the complete geodesic extending $\overline{x_n y_n}$ parameterized such that $\sigma_n(0)$ is the closest point to $\sigma(0)$. It is easy to see that for large values of $n$, $x_n$ and $y_n$ have to be on opposite sides of $\sigma_n(0)$. Therefore $d(\sigma_n(0), p)$ goes to infinity. Fix $R > 0$ and let $q_n$ be the point on $\overline{p\sigma_n(0)}$ such that $d(p, q_n) = R$. By passing to a subsequence we assume that $q_n$ converges to $q$ and that the rays $\overrightarrow{q_n x_n}$ and $\overrightarrow{q_n y_n}$ converge respectively to two rays starting at $q$ and asymptotic to $\overrightarrow{p\xi_x}$ and $\overrightarrow{p\xi_y}$. In fact these two rays form a complete geodesic, see [Bal82] for details. Therefore $\sigma$ is contained in a flat strip of width $R$. Since $R$ is arbitrary number, the geodesic $\sigma$ is contained in flat strips of arbitrary width. By a standard compactness argument, it is easy to show that a subsequence of the flat strips will converge to a flat half plane, which implies that $Td(\xi_x, \xi_y) = \pi$, and this is a contradiction. $\square$

For any convex subset $K$ of a Hadamard manifold $X$, we denote by $\pi_K$ the projection map from $X$ to $K$.

**Proposition 3.4.** *There exists a constant $\delta > 0$ which depends only on the analytic 4-manifold $X$, such that for any maximal higher rank submanifold $F$, any $q \in F$, and any $p \notin F$, $d(\pi_F(p), \overline{pq}) \leq \delta$. In particular the triangle $\Delta(p, q, \pi_F(p))$ is $\delta$-thin.*

*Proof.* Modulo the action of the fundamental group $\Gamma$, there are only a finite number of maximal higher rank submanifolds. If the proposition is false, then there exists a maximal higher rank submanifold $F$ and a sequence of triangles $\Delta(p_n, q_n, \pi_F(p_n))$ such that $d(\pi_F(p_n), \overline{p_n q_n}) \geq n$. Since $F$ is closed, we could assume that $\pi_F(p_n)$ is contained in a compact subset of $F$. By passing to a subsequence we assume that $p_n$, $q_n$, and $\pi_F(p_n)$ converge respectively to $\xi_p \in \partial_\infty X$, $\xi_q \in \partial_\infty F$, and $c \in F$. By Proposition 3.3, $\xi_p \in \partial_\infty F$, which is a contradiction. This finishes the proof. $\square$

**Lemma 3.5.** *There exists a constant $D_1 > 0$, which depends only on the space $X$ such that if $F$ is a maximal higher rank submanifold and $\overline{pq}$ is a geodesic segment, then either $d(\pi_F(p), \pi_F(q)) \leq D_1$ or $d(\overline{pq}, F) \leq D_1$.*

*Proof.* Assume that the statement is false, then there is a maximal higher rank submanifold $F$ such that for every $n \in \mathbb{N}$ there exists a geodesic segment $\overline{p_n q_n}$ such that $d(\pi_F(p_n), \pi_F(q_n)) \geq n$ and $d(\overline{p_n q_n}, F) \geq n$. Since $F$ is closed, we could assume that $\pi_F(p_n)$ is contained in a compact subset of $F$. By passing to a subsequence we assume that $\pi_F(p_n)$, $\pi_F(q_n)$, $p_n$, and $q_n$ converge respectively to $c \in F$, $\eta \in \partial_\infty F$, and $\xi_p, \xi_q \in \partial_\infty X$. It is not hard to see that since $F$ is a totally geodesic submanifold, $\xi_p \notin \partial F$. Since $d(c, \overline{\pi_F(q_n)q_n}) \geq n - 1$ for large values of $n$ and $d(c, \overline{p_n q_n}) \geq n$, using Proposition 3.3 we see that $\xi_p$ is path connected to $\eta \in \partial_\infty F$, and therefore $\xi_p \in \partial F$ which is a contradiction. This finishes the proof of the lemma. $\square$

**Corollary 3.6.** *There exists a constant $D_2 > 0$ such that if $F$ is a maximal higher rank submanifold and $\overline{pq}$ is a geodesic segment such $d(\pi_F(p), \pi_F(q)) > D_2$ then there are two points $p', q' \in \overline{pq}$ such that $d(p', \pi_F(p)) \leq D_2$ and $d(q', \pi_F(q)) \leq D_2$. In particular $\overline{pq}$ runs within $D_2$ distance from the path consisting of the geodesic segments $\overline{p\pi_F(p)}$, $\overline{\pi_F(p)\pi_F(q)}$, and $\overline{\pi_F(q)q}$.*

*Proof.* Let $D_2 = D_1 + \delta$, where $D_1$ is the constant in Lemma 3.5 and $\delta$ is the constant in Proposition 3.4. Let $r \in \overline{pq}$ be a point which is closest to $F$. By Lemma 3.5, $d(r, \pi_F(r)) \leq D_1$ and therefore the two geodesic segments $\overline{pr}$ and $\overline{p\pi_F(r)}$ are at most $D_1$ apart. By Proposition 3.4, there exists a point on $\overline{p\pi_F(r)}$ which is $\delta$ close to $\pi_F(p)$. Therefore there is a point on $\overline{pr}$ which is $D_2$ close to $\pi_F(p)$. Similarly there exists a point on $\overline{rq}$ which is $D_2$ close to $\pi_F(q)$ and the corollary follows. □

Now we start the proof of the main result in this section.

*Proof of Theorem 3.2.* Assume that the statement is false, then there exists a sequence of triangles $\Delta(p_n, q_n, r_n)$ which are not $n$-thin relative to any maximal higher rank submanifold. In particular they are not $n$-thin triangles. Therefore there exists a point $c_n$ on $\overline{p_n q_n}$ such that $d(c_n, \overline{p_n r_n}) \geq n$ and $d(c_n, \overline{q_n r_n}) \geq n$. Since $\Gamma$ acts cocompactly on $X$, we assume that $c_n$ are contained in a compact subset. By passing to a subsequence we could assume that $c_n$ converges to a point $c \in X$. Notice that $p_n$, $q_n$, and $r_n$ diverge to infinity. By passing further to a subsequence we could assume that $p_n$, $q_n$, and $r_n$ converge respectively to $\xi_p$, $\xi_q$, and $\xi_r$ in $\partial_\infty X$. By Proposition 3.3, $Td(\xi_p, \xi_r) \leq \pi$ and $Td(\xi_r, \xi_q) \leq \pi$ and therefore they belong to a connected component of the Tits boundary. Clearly $\xi_p \neq \xi_q$, therefore there exists a maximal higher rank submanifold $F$ such that $\xi_p, \xi_q, \xi_r \in \partial_T F$.

We claim that $c$ belongs to $F$. If not, then the geodesic segments $\overline{p_n q_n}$ would converge to a geodesic $\sigma$ passing through $c$ and parallel to $F$. Let $\sigma'$ be a geodesic in $F$ which is parallel to $\sigma$. By the analyticity of $X$, $\sigma$ and $\sigma'$ are contained in a 2-flat. If $F \notin \mathcal{W}$, i.e. isometric to $\mathbb{R}^2$ or $\mathbb{R}^3$ then the parallel set $P_{\sigma'}$ of $\sigma'$ properly contains $F$ which contradicts the maximality of $F$. If $F \in \mathcal{W}$, then $\sigma'$ can not be a singular geodesics in $F$, otherwise the parallel set $P_{\sigma'}$ is 4-dimensional and $X$ would be reducible. So, we assume that $\sigma'$ is not a singular geodesic in $F$. In this case $P_\sigma$ is 3-dimensional and therefore has to be of the form $Q \times \mathbb{R}$ and therefore in $\mathcal{W}$. But by the assumption on $X$ that can not happen. Otherwise $\partial_T X$ would have a nonstandard component.

Let $p'_n = \pi_F(p_n)$, $q'_n = \pi_F(q_n)$, and $r'_n = \pi_F(r_n)$. Without loss of generality assume that $d(q'_n, r'_n) \geq d(p'_n, r'_n)$, we might need to pass to a subsequence to guarantee that for every $n$. It is easy to see that $d(c, p'_n)$ and $d(c, q'_n)$ go to infinity since otherwise $p_n$ respectively $q_n$ would not converge to a point in $\partial_T F$. This implies that $d(p'_n, q'_n)$ goes to infinity, and therefore $d(q'_n, r'_n)$ goes to infinity. We need to show that $d(r'_n, p'_n)$ goes to infinity as well. If not then by passing to a subsequence we could assume that $d(p'_n, r'_n) \leq C$, for some constant $C$.

By Corollary 3.6, there exist four points $t_n, t'_n \in \overline{p_n q_n}$ and $s_n, s'_n \in \overline{q_n r_n}$ such that $d(t_n, p'_n)$, $d(t'_n, q'_n)$, $d(s_n, r'_n)$, and $d(s'_n, q'_n)$ are smaller than or equal $D_2$. Therefore $d(t_n, s_n) \leq 2D_2 + C$. Therefore the two geodesic segments $\overline{q_n s_n}$ and $\overline{q_n t_n}$ are at most $2D_2 + C$ apart. But this contradicts that $d(c, \overline{q_n r_n})$ goes to infinity. Therefore $d(p'_n, r'_n)$ goes to infinity as well. Again by Corollary 3.6, there are two point $l_n, l'_n \in \overline{p_n r_n}$ such that $d(l_n, p'_n)$ and $d(l'_n, r'_n)$ are smaller than or equal $D_2$. Now it is easy to see that all the triangles $\Delta(p_n, q_n, r_n)$ are $D_2$-thin relative to $F$, contradicting the choice of these triangles. This is finishes the proof of the theorem. □

*Remark* 3.7. The proof of Theorem 3.2 shows that the "fat" part of any triangle is close to a maximal higher rank submanifold.

## 4. The Asymptotic cone of $X$

In this section we analyze the asymptotic cone of the Hadamard manifold $X$. As in Section 3, we still assume that $X$ is the universal cover of a closed irreducible 4-dimensional nonpositively curved real analytic manifold without nonstandard components in $\partial_T X$.

We now recall the definition of the asymptotic cone of a metric space $(X, d)$. Fix a non-principle ultrafilter $\omega$ on $\mathbb{N}$, a sequence of base points $x_n \in X$, and a sequence of rescaling factors $\lambda_n$ such that $\omega\text{-}\lim \lambda_n = \infty$. The based ultralimit of $((X, \lambda_n^{-1} \cdot d), x_n)$ decomposes generically into uncountable number of components. The asymptotic cone is defined to be the component of the ultralimit containing the base point. We refer the reader to [KL95, KL97b] for further details.

We recall two well known facts about the asymptotic cone which we will use. For any non-principle ultrafilter $\omega$, $\text{Cone}_\omega(\mathbb{R}^n) = \mathbb{R}^n$. If $X$ is a quasi-homogeneous Gromov hyperbolic space with uncountable number of ideal boundary points, then $\text{Cone}_\omega(X)$ is an $\mathbb{R}$-tree with uncountable branching.

If $X$ is a Hadamard manifold, then $\text{Cone}_\omega(X)$ is a Hadamard space. Any sequence of flats in $X$ gives rise to a flat in the ultralimit of $((X, \lambda_n^{-1} \cdot d), x_n)$. If the distance between that flat and the base point in the asymptotic cone is finite, then the flat is in $\text{Cone}_\omega(X)$. The goal of this section is to show that these are the only flats which appear in $\text{Cone}_\omega(X)$.

Our analysis is similar to the analysis done by Kapovich and Leeb in [KL95]. And we will use some of their results. We will often mention the results without proof and refer the reader to their paper for the proofs.

Given a maximal higher rank submanifold $F \in \mathcal{V}$, if $F \approx \mathbb{R}^2$, $\mathbb{R}^3$, or $Q \times \mathbb{R}$, the sequence $F_n = F$ has a limit $\mathbb{R}^2$, $\mathbb{R}^3$, or $T \times \mathbb{R}$, where $T$ is an $\mathbb{R}$-tree, in $\text{Cone}_\omega(X)$.

The group $\Gamma^* = \prod_\mathbb{N} \Gamma$ acts by isometries on $\text{Cone}_\omega(X)$. Recall that up to the action of $\Gamma$ there is only a finite number of maximal higher rank submanifolds, we denote them by $K_1, \ldots, K_r$. Any sequence $F_n$ of maximal higher rank submanifolds in $\mathcal{V}$ gives rise to a partition $A_1 \sqcup \cdots \sqcup A_r$ of $\mathbb{N}$ as follows: $n \in A_s$ if $F_n = \phi(K_s)$ for some isometry $\phi \in \Gamma$. By basic properties of ultrafilters, there is exactly one subset $A_s \in \omega$. Therefore the sequence $F_n$ and the subsequence corresponding to $A_s$ give rise to the same limit in $\text{Cone}_\omega(X)$. Therefore without loss of generality we could assume that for a fixed $K_s$, $F_n = \phi_n(K_s)$, where $\phi_n \in \Gamma$. Since $\Gamma^*$ acts by isometries on $\text{Cone}_\omega(X)$, it is easy to see that the limit of $F_n$ is $\phi^*(\omega\text{-}\lim K_s)$ where $\phi^* = (\phi_1, \phi_2, \ldots) \in \Gamma^*$.

We denote by $\mathcal{F}$ the limits in $\text{Cone}_\omega(X)$ of all sequences of maximal higher rank submanifolds of $X$. The above discussion shows the following,

**Proposition 4.1.** *Every element in $\mathcal{F}$ is isometric to $\mathbb{R}^2$, $\mathbb{R}^3$, or $T \times \mathbb{R}$, where $T$ is an $\mathbb{R}$-tree.*

**Definition 4.2.** Let $X$ be a CAT(0) space. A triangle $\Delta(p_1, p_2, p_3)$ in $X$ is called an open triangle if $p_1$, $p_2$, $p_3$ are different points, and $\overline{p_i p_j} \cap \overline{p_i p_k} = \{p_i\}$, where $i, j, k \in \{1, 2, 3\}$ and $i \neq j \neq k$.

The following lemma is a rephrase of Proposition 4.3 in [KL95].

**Lemma 4.3.** *The asymptotic cone $\text{Cone}_\omega(X)$ satisfies the following properties:*

*1. Every open triangle is contained in some $F \in \mathcal{F}$.*

*2. Any two different elements $F, F' \in \mathcal{F}$ have at most one point in common.*

*Proof.* The proof of the first part is identical to the proof of the first part in Proposition 4.3 in [KL95], and the proof carries over to our sitting.

Now we give the proof of the second part. Notice that every $F \in \mathcal{F}$ is a convex subset of $\text{Cone}_\omega(X)$. Assume that $F$ and $F'$ intersect in two different points $x$ and $y$. Therefore $\overline{xy} \subset F \cap F'$. Let $F_n$ be a sequence of maximal higher rank submanifolds such that $\omega\text{-}\lim F_n = F$. Let $z \in F'$ be any point such that the triangle $\Delta(x, y, z)$ is an open triangle. The goal is to show that $z \in F$. Choose a sequence of triangles $\Delta_n = \Delta(x_n, y_n, z_n)$ in $X$ such that they converge to $\Delta(x, y, z)$ in the asymptotic cone. We can choose $x_n, y_n \in F_n$. Let $z'_n = \pi_{F_n}(z_n)$ be the projection of $z_n$ to $F_n$.

We claim that the two sequences $(z_n)_{n \in \mathbb{N}}$ and $(z'_n)_{n \in \mathbb{N}}$ represent the same point $z$ in $\text{Cone}_\omega(X)$, which implies that $z \in \omega\text{-}\lim F_n = F$. To see that, assume that $(z'_n)_{n\in\mathbb{N}}$ converges to $z' \in \mathrm{Cone}_\omega(X)$ and $z \neq z'$. Using Proposition 3.4, for every $n \in \mathbb{N}$, we can choose two points $s_n \in \overline{z_n x_n}$ and $t_n \in \overline{z_n y_n}$ within distance at most $\delta$ from $z'_n$, where $\delta$ is the constant given in Proposition 3.4. Clearly, $(z'_n)_{n\in\mathbb{N}}$, $(s_n)_{n\in\mathbb{N}}$, and $(t_n)_{n\in\mathbb{N}}$ represent the same point $z' \in \mathrm{Cone}_\omega(X)$ which implies that $\{z\} \subsetneqq \overline{zz'} \subset \overline{zx} \cap \overline{zy}$. This is a contradiction since we assumed that $\Delta(x, y, z)$ is an open triangle.

Notice that the set of points in $F'$ where the triangle $\Delta(x, y, z)$ is open is an open dense subset of $F'$, for the three different possibilities of $F'$ given by Proposition 4.1. By a continuity argument we see that $F' \subset F$. Similarly we show that $F \subset F'$ which finishes the proof. $\square$

*Remark* 4.4. The proof of Lemma 4.3 rules out the possibility that the limit of a sequence of maximal higher rank submanifolds which are isometric to $\mathbb{R}^2$ is contained in the limit of a sequence of maximal higher rank submanifolds of higher dimensions.

For every element $F \in \mathcal{F}$, we denote by $\pi_F \colon \mathrm{Cone}_\omega(X) \to F$ the projection map which is well defined and distance non-increasing since $F$ is a convex subset of $\mathrm{Cone}_\omega(X)$ which is CAT(0) space. Now we state Lemma 4.4 and Lemma 4.5 from [KL95] in our new sitting where 2-flats are replaced by maximal higher rank submanifolds.

**Lemma 4.5.** *Let $\gamma \colon I \to \mathrm{Cone}_\omega(X)$ be a curve in the complement of $F$. Then $\pi_F \circ \gamma$ is constant.*

**Lemma 4.6.** *Every embedded closed curve $\gamma \subset \mathrm{Cone}_\omega(X)$ is contained in some element $F \in \mathcal{F}$.*

The proofs of these two lemmas in [KL95] carry over to our new setting. The only ingredient used there was that the limits of 2-flats in the asymptotic cones of the universal covers of certain Haken manifolds, satisfy the two properties in Lemma 4.3. The proofs carry over to our new sitting after replacing 2-flats by maximal higher rank submanifolds and using Lemma 4.3.

Lemma 4.3 easily implies that the only flats in $\mathrm{Cone}_\omega(X)$ are the ones which are subsets of elements of $\mathcal{F}$. While Lemma 4.6 shows that every embedded disk of dimension at least 2 is contained in an element of $\mathcal{F}$.

## 5. Proof of the main theorem

In this section we give the proof of Theorem 1.1. We assume that $\partial_T X_1$ has a nonstandard component and $\partial_T X_2$ does not have any nonstandard components. If $f\colon X_1 \to X_2$ is a quasi-isometry, then it induces a map $\text{Cone}_\omega(f)\colon \text{Cone}_\omega(X_1) \to \text{Cone}_\omega(X_2)$ which is bi-Lipschitz homeomorphism. Since $\partial_T X_1$ contains a nonstandard component, there exist two maximal higher rank submanifolds $H_1 = Q_1 \times \mathbb{R}$ and $H_2 = Q_2 \times \mathbb{R}$ which intersect in a 2-flat, which we denote by $K$. In the asymptotic cone $Q_1 \times \mathbb{R}$ and $Q_2 \times \mathbb{R}$ give rise to two convex subsets $W_1 = T_1 \times \mathbb{R}$ and $W_2 = T_2 \times \mathbb{R}$, where $T_1$ and $T_2$ are $\mathbb{R}$-trees. The $\mathbb{R}$-factors of $W_1$ and $W_2$, which we call the singular directions are different. The 2-flat $K$ gives rise to a 2-flat, which we denote by $L$, in the asymptotic cone $\text{Cone}_\omega X_1$. Clearly $W_1 \cap W_2 \supseteq L$. The goal of the following proposition is to show that the intersection of $W_1$ and $W_2$ is precisely $L$.

**Proposition 5.1.** *Using the above notation, $W_1 \cap W_2 = L$.*

*Proof.* Let $z \in W_1 \cap W_2$. Choose two sequences $(x_n)_{n \in \mathbb{N}} \subset H_1$ and $(y_n)_{n \in \mathbb{N}} \subset H_2$ which represent the point $z$ in $\text{Cone}_\omega(X_1)$. This implies that $\omega\text{-}\lim d(x_n, y_n)/\lambda_n = 0$. Since $H_1$ and $H_2$ are orthogonal to each other at the intersection (see Lemma 3.4 in [Sch89b]), therefore $\pi_{H_1}(y_n) \in K = H_1 \cap H_2$. Since $d(\pi_{H_1}(y_n), y_n) \leq d(x_n, y_n)$, it easy to see that $(x_n)_{n \in \mathbb{N}}$, $(y_n)_{n \in \mathbb{N}}$, and $(\pi_{H_1}(y_n))_{n \in \mathbb{N}}$ represent the same point in $\text{Cone}_\omega(X_1)$. This shows that $z \in L$, which finishes the proof. $\square$

Since $\partial_T X_2$ does not contain any nonstandard components, Section 4 describe all the flats in $\text{Cone}_\omega(X_2)$. We are finally ready to start the proof of the main theorem.

*Proof of Theorem 1.1.* Assume that $f\colon X_1 \to X_2$ is a quasi-isometry. The induced map $\text{Cone}_\omega(f)\colon \text{Cone}_\omega(X_1) \to \text{Cone}_\omega(X_2)$ is a bi-Lipschitz homeomorphism. Assume that $T_i \times \mathbb{R} \subset \text{Cone}_\omega(X_1)$, for $i = 1, 2$, where $T_i$ is an $\mathbb{R}$-tree such that $T_1 \times \mathbb{R} \cap T_2 \times \mathbb{R} = \mathbb{R}^2$ as mentioned above. By Lemma 4.6 the image of each flat in $T_i \times \mathbb{R}$, for $i = 1, 2$, has to be contained inside some element $F \in \mathcal{F}$. For any two 2-flats in $T_i \times \mathbb{R}$, there exists a third flat which intersects both of them in half planes. That shows that the images of $T_i \times \mathbb{R}$, for $i = 1, 2$, have to be contained in some $F \in \mathcal{F}$. We need to show this is not possible for the three different types of $F$, namely $\mathbb{R}^2$, $\mathbb{R}^3$, and $T \times \mathbb{R}$.

*Case 1: $F = \mathbb{R}^2$.* This case is easy, a bi-Lipschitz embedding of any 2-flat in $T_1 \times \mathbb{R}$ has to be onto $F$, contradicting that $\text{Cone}_\omega(f)$ is a bijection.

*Case 2: $F = \mathbb{R}^3$.* Since $T_i = \text{Cone}_\omega(Q_i)$ where $Q_i$ are visibility 2-manifolds (with cocompact action) and therefore hyperbolic, then $T_i$ branches uncountably many times at each point. Now it is not hard to see that an $\mathbb{R}$-tree which branches uncountably many times at each point can not be bi-Lipschitz embedded in $\mathbb{R}^3$. One easy way to see that is to fix a base point $p$ in the tree, and take a sequence of points $p_i$ such that $d(p, p_i) = 1$ and $d(p_i, p_j) = 2$, for $i \neq j$. But $\text{Cone}_\omega(f)(p_i)$ is contained in a compact subset of $\mathbb{R}^3$ and therefore by passing to a subsequence it will converge in $\mathbb{R}^3$, which contradicts that $\text{Cone}_\omega(f)$ is bi-Lipschitz homeomorphism.

*Case 3:* $F = T \times \mathbb{R}$. First we recall Lemma 2.14 in [KL95] which states that if $T$ is a metric tree then the image of any bi-Lipschitz embedding $\pi\colon \mathbb{R}^2 \longrightarrow T \times \mathbb{R}$, is a flat in $T \times \mathbb{R}$. Notice that any two different flats in $F = T \times \mathbb{R}$ either do not intersect, intersect in a line, a strip, or a half plane.

Denote by $\alpha_1$, $\alpha_2$ the unique geodesics in $T_1$, $T_2$ respectively which satisfies that $\alpha_1 \times \mathbb{R} = \alpha_2 \times \mathbb{R}$ is the unique 2-flat which is the intersection of $T_1 \times \mathbb{R}$ and $T_2 \times \mathbb{R}$. Parameterize $\alpha_1$ and $\alpha_2$ such $p = \alpha_1(0) = \alpha_2(0)$ is the unique intersection point. Choose a geodesic $\beta_1$ in $T_1$ such that, $\beta_1$ follows $\alpha_1$ until they reach the point $p$ and then branches at $p$. Choose a geodesic $\beta_2$ in $T_2$ such the intersection of $\alpha_2$ and $\beta_2$ is the unique point $p$. Consider the two flats $\beta_1 \times \mathbb{R} \subset T_1 \times \mathbb{R}$ and $\beta_2 \times \mathbb{R} \subset T_2 \times \mathbb{R}$. It is not hard to see that these two planes intersect in a half line. But by Lemma 2.14 in [KL95] mentioned above, the image of these two planes in $F$ are two flats. This is a contradiction since a half line is not bi-Lipschitz homeomorphic to either a line, a strip or a half plane. $\square$

DEPARTMENT OF MATHEMATICS, UNIVERSITY OF PENNSYLVANIA, PHILADELPHIA, PA 19104-6395, USA

*E-mail address*: mhindawi@math.upenn.edu